\documentclass[12pt,a4paper]{amsart}
\usepackage{amssymb}

\title[Harmonic measure and uniform densities]
{Harmonic measure and uniform densities}
\author{Joaquim Ortega-Cerd\`a}
\address{Dept.\ Matem\`atica Aplicada i An\`alisi, Universitat de Barcelona,
Gran Via 585, 08071 Barcelona, Spain}
\email{jortega@ub.edu}
\thanks{The authors are supported by the European Commission Research Training
Network HPRN-CT-2000-00116. The first author is supported by the
DGICYT grant: BFM2002-04072-C02-01 and by the CIRIT grant:
2001SGR00172. The second author is partly supported by a grant
from the Research Council of Norway. }
\author{Kristian Seip}
\address{Dept.\ of Mathematical Sciences, Norwegian University of Science and
Technology, N--7491 Trondheim, Norway} \email{seip@math.ntnu.no}
\keywords{Harmonic measure, lower and upper uniform densities}
\subjclass{30C85, 31C05}
\date{\today}

\newcommand{\D}{\mathbb D}
\newcommand{\C}{\mathbb C}

\newtheorem{theorem}{Theorem}
\newtheorem*{theoremA}{Theorem A}
\newtheorem*{theoremB}{Theorem B}
\newtheorem{lemma}{Lemma}

\begin{document}
\begin{abstract}
We study two problems concerning harmonic measure on ``champagne
subdomains" of the unit disk $\D$. These domains are obtained by
removing from $\D$ little disks around sequences of points with a
uniform distribution with respect to the pseudohyperbolic metric
of $\D$. We find (I) a necessary and sufficient condition on the
decay of the radii of the little disks for the exterior boundary
to have positive harmonic measure, and (II) describe sampling and
interpolating sequences for Bergman spaces in terms of the
harmonic measure on such ``champagne subdomains''.
\end{abstract}
\maketitle
\section{Introduction}

This paper presents two theorems concerning harmonic measure on
``champagne subdomains'' of the unit disk. Our first result
solves a problem posed in a recent paper by Akeroyd
\cite{Akeroyd02}, while our second result gives a Bergman space
counterpart of a result of Garnett, Gehring, and Jones
\cite{GGJ83} for interpolation by bounded analytic functions.

The setting is as follows. Let $\Lambda$ be a sequence of distinct
points in the open unit disk $\D=\{z\in\C:\ |z|<1\}$, and define
\[
\rho(z,\zeta)=\left|\frac{z-\zeta}{1-\zeta\bar z}\right|,
\]
which is the pseudohyperbolic distance between $z,\zeta\in \D$.
For $z\in\D$ and $0<r<1$ we set
\[
D(z,r)=\{\zeta\in\D:\ \rho(z,\zeta)\le r\}.
\]
We say that
$\Lambda$ is a \emph{uniformly dense sequence} if
\begin{itemize}
\item[(i)] $\Lambda$ is separated , i.e., $\inf_{\lambda\ne\lambda'}
\rho(\lambda,\lambda')>0$, $\lambda,\lambda'\in\Lambda$.
\item[(ii)] There exists an $r<1$ such that $\D=\bigcup_{\lambda\in\Lambda}
D(\lambda,r)$.
\end{itemize}
Uniformly dense sequences appear naturally in the study of Bergman
spaces, as sampling or interpolating sequences (see below).

We are interested in studying harmonic measure on ``champagne
subdomains'' with ``bubbles" around the points
$\lambda\in\Lambda$, i.e., we consider infinitely connected
domains of the form
\begin{equation} \label{domain}
\Omega=\D\setminus \bigcup_{\lambda\in\Lambda} D(\lambda,r_\lambda),
\end{equation}
where the $r_\lambda$ ($0<r_\lambda<1$) are such that the closed
disks $D(\lambda,r_\lambda)$ are pairwise disjoint. The notation
$\omega(z, A, \Omega)$ stands for the value at $z$ of the harmonic
measure on $\Omega$ of a set $A\in
\partial\Omega$. (See below for a formal definition.) We will be concerned with the
following two problems:
\begin{itemize}
\item[(I)] Find a necessary and sufficient condition on the decay
of $r_\lambda$ for the exterior boundary $\partial\D$ to have
positive harmonic measure, i.e., for $\omega(z,\partial\D,
\Omega)>0$ to hold when $z\in \Omega$. \item[(II)] Characterize
sampling and interpolating sequences in terms of harmonic measure
on such ``champagne subdomains''.
\end{itemize}

We were led to Problem (I) by a question in \cite{Akeroyd02}. The
issue was whether for any uniformly dense sequence $\Lambda$
one may pick the $r_\lambda$ such that the
exterior boundary $\partial\D$ has zero harmonic measure and
\begin{equation}\label{rectify}
\sum_{\lambda\in\Lambda} \text{length}(\partial D(\lambda,r_\lambda)) <\infty.
\end{equation}
Our solution shows that there is ample room for such constructions.

When dealing with Problem (I), we assume that
$r_\lambda=\varphi(|\lambda|)$ with $\varphi$ a nonincreasing
function bounded by some constant less than $1$. In Theorem~\ref{dichotomy}
below, we arrive at the  following  necessary and sufficient
condition for $\partial\D$ to have positive harmonic measure:
\begin{equation}\label{logint}
\int_0^1 \frac{dt}{(1-t)\log(1/\varphi(t))} <\infty.
\end{equation}
Thus, for example, $\varphi(t)=c(1-r)^\gamma$ for arbitrary $\gamma>1$  yields
zero harmonic measure of $\partial\D$ as well as \eqref{rectify}. Note that the
result is independent of the particular choice of $\Lambda$.

Both problems (I) and (II) can be seen as originating from a paper
by Garnett, Gehring, and Jones \cite{GGJ83}, dealing with similar
considerations when $\Lambda$ satisfies the Blaschke condition
\begin{equation}\label{blaschke0}
\sum_{\lambda\in\Lambda} (1-|\lambda|)<\infty. \end{equation} On
the one hand, our work contrasts the following elementary fact:
Suppose $\Lambda$ is merely separated and $r_\lambda=r<1$ for all
$\lambda\in\Lambda$. In this case, if $\Lambda$ satisfies
\eqref{blaschke0}, then the exterior boundary of $\Omega$ has
positive harmonic measure. This can be seen as a statement about
the sparsity of Blaschke sequences, while our condition
\eqref{logint} reflects the density of uniformly dense sequences.
On the other hand, our work parallels \cite{GGJ83} in that we give
a description of  interpolating sequences for Bergman spaces in
terms of the harmonic measure (see Theorem~\ref{denshar} below), in a similar
way as is obtained for classical interpolating sequences in
\cite{GGJ83}.

It is interesting to note that the integral condition \eqref{logint} has
appeared before, in a different context. In \cite{LyuSei97b}, a sequence
$Z$ of distinct points in $\D$ is said to be a separated
non-Blaschke sequence if $Z$ is separated and
\[
\sum_{z\in Z} (1-|z|)=\infty.
\]
Also, $\varphi$ (again a nonincreasing function bounded by some constant less
than $1$)  is an \emph{essential minorant for} $H^\infty$ if the inequality
\[
|f(z)|\le \varphi(|z|) \ \ \text{for all} \ z\in Z,
\]
$f$ a bounded analytic function and $Z$  some separated
non-Blaschke sequence, implies that $f\equiv 0$. The theorem
proved in \cite{LyuSei97b} says that $\varphi$ is an essential
minorant for $H^\infty$ if and only if \eqref{logint} holds. We
are not able to give a direct proof of the link between essential
minorants and  our ``champagne subdomains'' whose exterior
boundaries have positive harmonic measure,  but we will offer a
heuristic argument for the connection.

The following notation will be used repeatedly below:
We write  $A \lesssim B$ to signify that $A\le C B$ for some
constant $C>0$, independent of whatever arguments are involved. If
both $A\lesssim B$ and $B\lesssim A$, then we write $A\simeq B$.

\section{Positive harmonic measure of the exterior boundary}

We begin by noting that the domains $\Omega$ defined by
\eqref{domain} are Dirichlet domains because they satisfy the
exterior cone condition. Thus every continuous function on
$\partial \Omega$ can be extended continuously to a harmonic
function in the interior of $\Omega$. The maximum principle shows
that the evaluation at any point $z\in \Omega$ of this extension
is a bounded linear functional on $\mathcal C(\partial \Omega)$
with norm less than $1$. By the Riesz representation theorem,
there is a probability measure $\omega_z$ supported on $\partial
\Omega$ such that the action of this functional on $f$ can be
represented as an integral against $\omega_z$. The measure
$\omega_z$ is called the harmonic measure of $\Omega$ at $z$, and
the harmonic measure of a set $A\subset
\partial \Omega$ is denoted $\omega(z,A,\Omega)$.

Note that the function $z\mapsto \omega(z,A,\Omega)$ is a
nonnegative  harmonic function on $\Omega$. The maximum principle
implies that this function is either identically $0$ or strictly
positive on $\Omega$.

There exist several equivalent definitions of harmonic measure.
One such definition is given in probabilistic terms: The harmonic
measure $\omega(z,A,\Omega)$ coincides with the probability that
a Brownian motion starting at the point $z$ exits the open set
$\Omega$ for the first time at one of the points in $A$.  We refer
to \cite{Bass95} for a proof of this fact and for some examples of
estimates of the harmonic measure using this probabilistic
interpretation.

Fix a uniformly dense sequence $\Lambda$, and assume that $\varphi$ is a
nonincreasing function on $(0,1)$ such that the closed disks
$D(\lambda,\varphi(|\lambda|))$, $\lambda\in\Lambda$, are pairwise disjoint.
Set
\[
\Omega(\Lambda,\varphi)=\D\setminus\bigcup_{\lambda\in\Lambda}
D(\lambda,\varphi(|\lambda|)),
\]
and assume for convenience that $0\in \Omega(\Lambda,\varphi)$.

\begin{theorem}\label{dichotomy}
The exterior boundary of $\Omega(\Lambda,\varphi)$ has positive
harmonic measure, i.e.
\[
\omega(0,\partial\D,\Omega(\Lambda,\varphi))>0,
\]
if and only if
\begin{equation} \label{logint2}
\int_{0}^1 \frac{dt}{(1-t)\log(1/\varphi(t))} <\infty.
\end{equation}
\end{theorem}
Note that the condition may be written equivalently as
\begin{equation} \label{alternative}
\sum_{j=1}^\infty \frac{1}{\log(1/\varphi(1-K^{-j}))} < \infty
\end{equation} for some $K>1$.

Theorem~\ref{dichotomy} reflects the following dichotomy: Either the little
disks  are so small that the contribution to the harmonic measure from each of
them can be viewed as independent of the  contributions from the others, or the
disks are so large that their  contributions to the harmonic measure interact
in a profound way. The first case corresponds to positive harmonic measure of
the exterior boundary, the second case to zero harmonic measure of the exterior
boundary.

The proof of the sufficiency of \eqref{logint2} illuminates this point: We
begin by observing that we may safely disregard a finite number of points; thus
we may consider instead ($r<1$)
\[
\Lambda_r =\Lambda \cap \{z: \ |z|>r\}.
\]
It is immediate that
\begin{equation}\label{onehole}
\omega(0,\partial D(\zeta,s), \D \setminus D(\zeta,s))=
\frac{\log|\zeta|}{\log s}.
\end{equation}
Then
\[
1-\omega(0,\partial\D,\Omega(\Lambda_r,\varphi))
\le \sum_{|\lambda|\ge r} \frac{\log \frac{1}{|\lambda|}}
{\log \frac1{\varphi(|\lambda|)}}\lesssim\int_{r-(1-r)}^1
\frac{dt}{(1-t)\log 1/\varphi(t)},
\]
where the latter inequality follows from the fact that $\Lambda$ is a
separated sequence.
We are done because the integral can be made smaller than 1
by choosing $r$ sufficiently close to 1.

Before proving the necessity of \eqref{logint2}, we comment on the relation to
essential minorants. As explained in \cite{LyuSei97b}, the  result describing
essential minorants is really a statement about the size of the exceptional set
on which a positive superharmonic function $u(z)$ exceeds
$\log(1/\varphi(|z|))$. It
was proved in \cite{LyuSei97b} that  with $m$ denoting Lebesgue area measure on
$\D$, we have
\begin{equation}\label{except}
\int_{u(z)>\log(1/\varphi(|z|))} \frac{dm(z)}{1-|z|} <\infty
\end{equation}
for each positive superharmonic function $u$
if and only if \eqref{logint2} holds.
By Harnack's inequality, we may recast the integral in \eqref{except} as
a sum:
\begin{equation} \label{blaschke}
\sum_{u(\lambda)>\log(1/\varphi(|\lambda|))} (1-|\lambda|) < \infty.
\end{equation}
Now if the little disks can be seen as acting independently of each other,
then again by Harnack's inequality and the Riesz representation formula,
\[
u(z)\gtrsim \sum_{u(\lambda)>\log(1/\varphi(|\lambda|))} (1-|\lambda|),
\]
and it follows that $\varphi$ is an essential minorant.

We now turn to the necessity of \eqref{logint2}. We will estimate
the probability that a Brownian motion starting at $0$ and moving
in $\Omega(\Lambda,\varphi)$ will reach $\partial \D$. (We assume
the motion is stopped once the particle exits
$\Omega(\Lambda,\varphi)$.) Define
\[
C_j=\{z: |z|=1-K^{-j}\},
\]
$j=0,1,\ldots$ and $K$ is some large constant chosen such that for
every $z\in C_{j-1}$ there is a nearby point $\lambda_z\in \Lambda$ in
the annulus bounded by $C_{j-1}$ and $C_j$ such that
\[
\sup_j\sup_{z\in C_j}\rho(z,\lambda_z)<1.
\]
%We may assume that the
%sequence $\Lambda$ is bounded away from the circles $C_j$.
%(Alternatively, we can make small loops around the points from
%$\Lambda$ and still call the resulting curves $C_j$.)
Let $P_j$ denote the probability that our Brownian motion hits
$C_j$. If $Q_j$ denotes the supremum of the probabilities that a
Brownian motion starting from some point at $C_{j-1}$ hits $C_j$,
then we get
\[
P_j\le Q_j P_{j-1},
\]
and so by induction
\[
P_n\le \prod_{j=1}^{n} Q_j.
\]
Thus it is necessary that $\prod_{j=1}^\infty Q_j >0$.
Equivalently, we have
\begin{equation}\label{Qprod}
\sum_{j=1}^\infty (1-Q_j) <\infty.
\end{equation}

Note that $1-Q_j$ is the infimum of the probabilities that a
Brownian motion starting from some point on $C_{j-1}$ hits a disk
$D(\lambda,\varphi(|\lambda|))$ before reaching $C_j$. For any
point on $C_{j-1}$ we may discard all disks except
$D_{\lambda_z}=D(\lambda_z,\varphi(|\lambda_z|))$ lying in the
annulus bounded by $C_{j-1}$ and $C_j$, because we are thus
diminishing the probability of hitting the disks. Therefore, if we
denote by $D_j$ the disk bounded by $C_j$, then
\[
1-Q_j \ge \inf_{z\in C_{j-1}} \omega (z,\partial D_{\lambda_z},D_j\setminus
D_{\lambda_z}).
\]
This harmonic measure can be estimated because $\sup\rho(z,\lambda_z)<1$:
\[
\omega (z,\partial D_{\lambda_z},D_j\setminus
D_{\lambda_z})\gtrsim \frac 1{\log (1/\varphi(1-K^{-j}))} \qquad \forall z\in
C_{j-1}.
\]
Combining this estimate with \eqref{Qprod}, we arrive at
\eqref{alternative}.

\section{Lower and upper uniform densities}\label{densitats}

In the previous section, the particular choice of uniformly dense
sequence $\Lambda$ was inessential. However, such sequences may
have different densities, and a natural question is whether these
densities can be captured in terms of harmonic measure. We will
now show how this can be done.

Let $\Lambda$ be a separated sequence. Following \cite{Seip93}, we
define the \emph{lower uniform density} of $\Lambda$ as
\[
D^{-}(\Lambda) = \liminf_{r\to 1} \inf_{z\in\D}
\frac{\sum_{\rho(\lambda,z)<r}(1-\rho(\lambda,z))}{\log
\frac1{1-r}}
\]
and \emph{the upper uniform density} of $\Lambda$ as
\[
D^{+}(\Lambda) = \limsup_{r\to 1} \sup_{z\in\D}
\frac{\sum_{\rho(\lambda,z)<r}(1-\rho(\lambda,z))}{\log
\frac1{1-r}}.
\]
Note that we always have $D^-(\Lambda)\le D^+(\Lambda)<\infty$,
and that $D^-(\Lambda)>0$ if and only if $\Lambda$ is a uniformly
dense sequence.

To see the significance of these densities, we cite the main
results of \cite{Seip93}. Let $A^{-\alpha}$ ($\alpha>0$) be the
space of analytic functions $f$ on $\D$ satisfying
\[
\| f \|_\alpha = \sup_{z\in \D} (1-|z|^2)^\alpha |f(z)|<\infty.
\]
We say that $\Lambda$ is a \emph{sampling sequence} for
$A^{-\alpha}$ if there is a positive constant $C$ such that
\[
\| f \|_\alpha\le C
\sup_{\lambda\in\Lambda}(1-|\lambda|^2)^\alpha|f(\lambda)|
\]
for every function $f\in A^{-\alpha}$. On the other hand, we say
that $\Lambda$ is an \emph{interpolating sequence} for
$A^{-\alpha}$ if the interpolation problem
\[
f(\lambda)=a_\lambda
\]
has a solution $f\in A^{-\alpha}$ whenever
$\{(1-|\lambda|^2)^\alpha a_\lambda\}$ is a  bounded sequence. In
\cite{Seip93}, it was proved that a separated sequence $\Lambda$
is a sampling sequence for $A^{-\alpha}$ if and only if
\[
D^{-}(\Lambda)>\alpha,
\]
and that $\Lambda$ is an interpolating sequence for $A^{-\alpha}$
if and only if
\[
D^{+}(\Lambda)<\alpha.
\]
These density conditions also describe similar sequences of
sampling and interpolation for weighted Bergman $L^p$ spaces.

Before stating our second theorem, we mention the result on which
it is modelled. A sequence $\Lambda$ of distinct points in $\D$ is
an interpolating sequence for $H^\infty$ if the interpolation
problem
\[
f(\lambda)=a_\lambda
\]
has a solution $f\in H^\infty$ whenever
$\{a_\lambda\}_{\lambda\in\Lambda}$ is a bounded sequence. We let
$\Lambda_\lambda$ be the sequence obtained from $\Lambda$ by
removing the one element $\lambda$. In \cite{GGJ83}, the following
re-interpretation of Carleson's theorem \cite{Carleson58} was
given:

\begin{theoremB}
A separated sequence $\Lambda$ is an interpolating sequence for
$H^\infty$ if and only if
\[
\inf_{\lambda\in\Lambda} \omega(\lambda, \partial \D,
\Omega(\Lambda_\lambda, c))>0
\]
for some $0<c<1$.
\end{theoremB}

To obtain a counterpart of this result, we define the following
densities. Set
\[
\Omega(z,r)=
\Omega(\Lambda;z,r)=\D\setminus \bigcup_{1/2<\rho(\lambda,z)<r}
D(\lambda,1-r),
\]
which is a finitely connected domain. We see that the uniform
pseudohyperbolic radius of the little disks tends to 0 as $r\to 1$.
This decay is tuned with the growth of $r$ in such a way that the
numbers
\[
D^{-}_h(\Lambda)=\liminf_{r\to 1^-}\inf_{z\in\D}
\log\frac{1}{\omega(z,\partial\D,\Omega(z,r))}
\]
and
\[
D^{+}_h(\Lambda)=\limsup_{r\to
1^-}\sup_{\lambda\in\Lambda}\log\frac{1}
{\omega(\lambda,\partial\D,\Omega(\lambda,r))}
\]
are positive when $\Lambda$ is uniformly dense. In fact, we have
the following precise characterization.
\begin{theorem}\label{denshar}
For a separated sequence $\Lambda$ in $\D$ we have
\[
D^-(\Lambda)=D^-_h(\Lambda) \ \ \text{and} \ \
D^+(\Lambda)=D^+_h(\Lambda).
\]
\end{theorem}

The proof of Theorem 2 combines probabilistic arguments with
certain precise function theoretic constructions, to be given in
the next section.

\section{Growth of analytic functions vanishing on $\Lambda$}

We will now see how $D^-(\Lambda)$ and $D^+(\Lambda)$ are related
to the growth of analytic functions vanishing on $\Lambda$. The
growth estimates to be established rely on a basic approximation
result for subharmonic functions.

We require some notation. If $f$ is analytic in $\D$, we denote by
$Z(f)$ its sequence of zeros. If $f$ has a zero of order $n$ at
$z$, then this is recorded by letting $z$ appear $n$ times in
$Z(f)$. On the other hand, we also think of $Z(f)$ as a subset of
the disk. In particular, when we say  that $Z(f)$ is separated, we
mean that $Z(f)$ consists of distinct points and that
\[
\inf_{z\neq z'}\rho(z,z')>0, \quad z,z'\in Z(f).
\]
We will rely on the following approximation result from
\cite{Seip95b}. 

\begin{theoremA}
Let $\Psi$ be subharmonic in $\D$ so that its Laplacian $\Delta
\Psi$ satisfies
\begin{equation} \label{nine}
\Delta \Psi (z)\simeq \frac{1}{(1-|z|^2)^2}
\end{equation}
for all $z\in \D$. Then there exists a function $g$ analytic in
$\D$, with $Z(g)$ a uniformly dense sequence, and
\[
|g(z)|\simeq \rho(z,Z(g)) e^{\Psi(z)}
\]
for all $z\in\D$.
\end{theoremA}

We deduce two lemmas from Theorem A.
\begin{lemma}\label{divisor}
Let $\Lambda$ be a uniformly dense sequence satisfying
$D^{-}(\Lambda)>\alpha>0$, and let $f$ be an analytic function on
$\D$ with $Z(f)=\Lambda$. Then there exists a uniformly dense
sequence $\Sigma$ and an analytic function $g$ on $\D$ with
$Z(g)=\Sigma$ such that
\[
\frac{|f(z)|}{|g(z)|}\simeq \frac{\rho(z,\Lambda)}{\rho(z,\Sigma)}
(1-|z|^2)^{-\alpha}
\]
for all $z\in \D$.
\end{lemma}
\begin{proof}
We proceed as in \cite{BerOrt95}. Set $u=\log|f|$,
\[
\xi_r(z)=(1-|z|)\chi_{D(0,r)}(z),
\]
and smooth $u$ by taking the invariant convolution
\[
u_r(w) = u\star \xi_r (w)= \frac 1{c_r}\int_{|z|<r}
u\left(\frac{z-w}{1-\overline{w}z}\right) (1-|z|)
\frac{dm(z)}{\pi(1-|z|^2)^2},
\]
where the constant $c_r$ is chosen in such a way that
$\|(1-|z|^2)^{-2}\xi_r\|_{L^1(\D)}=\pi c_r$. It follows (see
\cite{BerOrt95} for details) that
\begin{equation} \label{away}
|u_r(z)-u(z)|\le C_\varepsilon
\end{equation} if $\rho(z,\Lambda)>\varepsilon$,
and
\[
\frac{1}{4} (1-|z|^2)^2\Delta u_r(z) \ge \nu> \alpha
\]
for every $z$, provided $r$ is sufficiently close to $1$. The
latter inequality is a consequence of the assumption that
$D^{-}(\Lambda)>\alpha$. To see this, we use that the invariant
Laplacian $\tilde{\Delta}=(1-|z|^2)^2\Delta$ commutes with the
invariant convolution:
$$ \tilde{\Delta} u_r = (\tilde{\Delta} u)\star \xi_r. $$
This commutation holds because $\xi_r$ is a radial function.

We now consider the subharmonic function
\[
\Psi(z)= u_r(z)-\alpha\log\frac{1}{1-|z|^2},
\]
whose Laplacian satisfies \eqref{nine}. This means that Theorem A
applies. We obtain the desired estimates because $u$ and $u_r$ are
related by \eqref{away}; the estimates close to the points
$\lambda\in\Lambda$ follow from the maximum principle.
\end{proof}

Acting similarly as above, but considering instead
\[
\Psi(z)= \alpha\log\frac{1}{1-|z|^2} - u_r(z),
\]
we arrive at the following lemma.

\begin{lemma}\label{multiplicador}
Let $\Lambda$ be a separated sequence satisfying
$D^{+}(\Lambda)<\alpha$, and let $f$ be an analytic function on
$\D$ with $Z(f)=\Lambda$. Then there exists a uniformly dense
sequence  $\Sigma$ and an analytic function $g$ on $\D$ with
$Z(g)=\Sigma$ such that
\[
|f(z)g(z)|\simeq \rho(z,\Lambda)\rho(z,\Sigma) (1-|z|^2)^{-\alpha}
\]
for all $z\in \D$.
\end{lemma}

We note that the construction in \cite{Seip95b} can be adjusted so
that $\Lambda\cup \Sigma$ becomes a separated sequence in both
cases. However, we will not need this separation in what follows.

\section{Proof of Theorem 2}
Set $\alpha=D^-(\Lambda)$. We begin by showing that $\alpha\le
D^-_h(\Lambda)$. We prefer to give a general argument and define
\[
\Omega_\delta(z,r)=\D\setminus
\bigcup_{1/2<\rho(\lambda,z)<r}
D(\lambda,\delta(r)),
\]
where we only assume $\delta(r)\to 0$. Pick some small $\varepsilon>0$.  Let
$h=f/g$ be the function with zeros $\Lambda$ and poles $\Sigma$ given by
Lemma~\ref{divisor} such that
\[
|h(z)|\simeq (1-|z|)^{-\alpha+\varepsilon}
\]
far from $\Lambda$ and $\Sigma$; the constants involved here will  depend on
$\varepsilon$.

By conformal invariance, we may assume $z=0$. We will give a
probabilistic argument, estimating the probability that a Brownian
motion starting at $0$ and moving in $\Omega_\delta(0,r)$ will
reach $\partial D(0,r)$.  To this end, choose some function
$\eta(r)\to 0$, such that
\[
\log\frac{1}{\eta(r)}=o\left(\log\frac{1}{\delta(r)}\right),
\]
and
\[
\log\frac{1}{1-r}=n\log\frac{1}{\eta(r)}
\]
and define
\[
C_j(r)=\{z: |z|=1-\eta^{j}(r)\},
\]
$j=0,1,\ldots,n$. We may assume the sequence $\Sigma \cup \Lambda$
is bounded away from  the circles $C_j(r)$. (Alternatively, we can
make small loops around the points from $\Sigma$ and $\Lambda$.)
Let $P_j$ denote the probability that our Brownian motion hits
$C_j(r)$. If $Q_j$ denotes the supremum of the probabilities that
a Brownian motion starting from some point at $C_{j-1}(r)$ hits
$C_j(r)$, then we get
\[
P_j\le Q_j P_{j-1}
\]
and so by induction
\[
P_n\le \prod_{j=1}^{n} Q_j.
\]

To estimate $Q_j$, we disregard the points from $\Lambda$ on the inside of
$C_{j-1}(r)$. We may also disregard the points from $\Lambda$ close to
$C_j(r)$  (correspondingly, we divide out these zeros from $h$, but still call
the function $h$). This will increase the probability of hitting $C_{j}(r)$.
For  some constant $C$ (independent of $r$) the subharmonic function
\[
U_j(z)=\frac{\log\frac{1}{|h(z)|}+j(\alpha-\varepsilon)\log\frac{1}{\eta(r)}-C}
{\log\frac{1}{\delta(r)}+(\alpha-\varepsilon)\log\frac{1}{\eta(r)}}
\]
is bounded by $0$ on $C_j(r)$ and by $1$ on $\partial D(\lambda,\eta(r))$ for
$\lambda\in\Lambda$ between $C_{j-1}(r)$ and $C_j(r)$. Also, on $C_{j-1}(r)$ we
have
\[
 U_j(z)\ge
 \frac{(\alpha-\varepsilon)\log\frac{1}{\eta(r)}-2C}
 {\log\frac{1}{\delta(r)}+
        (\alpha-\varepsilon)\log\frac{1}{\eta(r)}}.
\]
It follows that
\[
Q_j\le \frac{\log\frac{1}{\delta(r)}+2C}{\log\frac{1}{\delta(r)}+
        (\alpha-\varepsilon)\log\frac{1}{\eta(r)}}.
\]
Thus
\[
\begin{split}
&\log \frac{1}{\omega(\lambda,\partial\D,\Omega(\lambda,r))}
\ge n\log \frac{\log\frac{1}{\delta(r)}+
    (\alpha-\varepsilon)\log\frac{1}{\eta(r)}}{\log\frac{1}{\delta(r)}+2C} \\
& = n\frac{(\alpha-\varepsilon)\log\frac{1}{\eta(r)}-2C}
          {\log\frac{1}{\delta(r)}+2C}
(1+o(1)) = (\alpha-\varepsilon) \frac{\log\frac{1}{1-r}}
                                     {\log\frac{1}{\delta(r)}}(1+o(1)).
\end{split}
\]

We now prove $\alpha \ge D^-_h(\Lambda)$. This time we cannot
allow $\delta(r)$ to decrease too slowly. Let $B_\zeta$ denote the
finite Blaschke product with zeros $\lambda\in\Lambda$ such that
$1/2<\rho(\zeta,\lambda)<r$. Define
\[
c(r)=\sup_z \log |B_z(z)|
\]
and pick $\zeta$ such that
\begin{equation}\label{bolet}
\log|B_\zeta(\zeta)| > c(r) -1.
\end{equation}
By conformal invariance, we may assume $\zeta=0$.
We set
\[
c(r)-1=-(\alpha-\varepsilon)\log\frac{1}{1-r}
\]
and note that by our definition of $\alpha$, $\varepsilon=\varepsilon(r)\to 0$
when $r\to 1$.

We introduce a function $\eta(r)$ as above and a similar
partition: Let $B_j$ be the Blaschke product with zeros
$\lambda\in\Lambda$ such that
\[
1-\eta^{j-1}(r)\le |\lambda| < 1-\eta^j(r).
\]
We define
\[
U_j(z)=\log \frac{1}{|B_j(z)|}.
\]
We now build a harmonic function which exceeds the harmonic
measure of the inner boundary. This function will be of the form
\[
U=\sum_{j=1}^n w_j U_j,
\]
with appropriate positive weights $w_j$ such that $U(z)\ge 1$ on
$\partial D(\lambda,\delta(r))$. To determine the $w_j$, we begin
by noting that
\[
U_n(z)\ge \log\frac{1}{\delta(r)}=: a
\]
on the boundary of the ``bubbles'' corresponding to the zeros of
$B_n$. Thus we set
\[
w_n=\frac{1}{a}.
\] Next observe that on the boundary of the ``bubbles" corresponding
to the zeros of $B_{n-1}(z)$, we get
\[
w_{n-1}U_{n-1}(z)+w_n U_n(z)\ge w_{n-1}a+ w_n
(\alpha-\xi)\log\frac{1}{\eta(r)}
\]
with $\xi=\xi(r)\to 0$ as $\eta(r)\to 0$. We set
\[
b:= (\alpha-\xi)\log\frac{1}{\eta(r)}
\]
and then
\[
w_{n-1}=\frac{a-b}{a^2}.
\]
Inductively, we get
\[
w_{n-j}=\frac{1}{a}\left(\frac{a-b}a\right)^j.
\]
To estimate $U(0)$, we argue as follows. The worst case is that
$U_n(0)$ is maximal because $w_n$ is the largest weight. Combining
our upper estimate \eqref{bolet}, i.e.,
\[
\sum_{j=1}^n U_j(0)\le (\alpha-\varepsilon) \log\frac{1}{1-r}
\]
with the lower estimates $U_j(0)\ge (\alpha-\xi)\log(1/\eta(r))$,
we get
\[
U_n(0)\le (\alpha-\xi)\log\frac{1}{\eta(r)}+
n(\xi-\varepsilon)\log\frac{1}{\eta(r)}.
\]
If $U_n(0)$ attains this upper bound, then
$U_j(0)=(\alpha-\xi)\log(1/\eta(r))=b$ for $1\le j <n$, and we
arrive at the estimate
\[
U(0)\le n(\xi-\varepsilon)\frac{\log\frac{1}{\eta(r)}}
{\log\frac{1}{\delta(r)}} + 1-\left(\frac{a-b}{a}\right)^n,
\]
and so
\[
\omega(0, \partial \D,\Omega_\delta(0,r))\ge
\left(\frac{a-b}{a}\right)^n -
(\xi-\varepsilon)\frac{\log\frac{1}{1-r}}{\log\frac{1}{\delta(r)}};
\]
We now require the second term to be ``small oh" of the first
term. This is certainly the case if $\delta(r)=1-r$. Thus
\[
\log\frac{1}{\omega(0,\partial\D,\Omega_\delta(0,r))}\le n
\frac{b}{a}(1+o(1))=\alpha \frac{\log\frac{1}{1-r}}
{\log\frac{1}{\delta(r)}} (1+o(1)),
\]
and we are done.

We next set $\alpha=D^+(\Lambda)$. The scheme is very similar. We
first show that $\alpha\ge D^+_h(\Lambda)$. Again
\[
\Omega_\delta(z,r)=\D\setminus
\bigcup_{\begin{subarray}{c}1/2<\rho(\lambda,z)<r\\ \lambda\in\Lambda
\end{subarray}}
D(\lambda,\delta(r)),
\]
where we only assume $\delta(r)\to 0$.
Pick some small $\varepsilon>0$.
Let $h=fg$ be the function  given by Lemma~\ref{multiplicador}
with zeros $\Lambda$
and $\Sigma$ such that
\[
|h(z)|\simeq (1-|z|)^{-\alpha-\varepsilon}
\]
far from $\Lambda$ and $\Sigma$; the constants involved here will
depend on $\varepsilon$.

By conformal invariance, we may assume $\lambda=0$. We estimate
again the probability that a Brownian motion starting at $0$ and
moving in $\Omega_\delta(0,r)$ will reach $\partial D(0,r)$. To
this end, choose some function $\eta(r)\to 0$, such that
\[
\log\frac{1}{\eta(r)}=o\left(\log\frac{1}{\delta(r)}\right),
\]
and
\[
\log\frac{1}{1-r}=n\log\frac{1}{\eta(r)}
\]
and define
\[
C_j(r)=\{z: |z|=1-\eta^{j}(r)\}
\]
$j=0,1,\ldots,n$. We may assume the sequence $\Lambda\cup\Sigma$
is bounded away from the circles $C_j(r)$. (Alternatively, we may
make small loops around the points.) Let $P_j$ denote the
probability that our Brownian motion hits $C_j(r)$. If $R_j$
denotes the infimum of the probabilities that a Brownian motion
starting from some point at $C_{j-1}(r)$ hits $C_j(r)$, then we
get
\[
P_j\ge R_j P_{j-1}
\]
and so by induction
\[
P_n\ge \prod_{j=1}^{n} R_j.
\]

We estimate $R_j$. For some constant $C$ (independent of $r$) the
superharmonic function
\[
U_j(z)=\frac{\log\frac{1}{|h(z)|}+j(\alpha-\varepsilon)\log\frac{1}{\eta(r)}+C}
{\log\frac{1}{\delta(r)}+(\alpha+\varepsilon)\log\frac{1}{\eta(r)}}
\]
is bounded from below by $0$ on $C_j(r)$ and by $1$ on
$\partial D(\lambda,\eta(r))$
for $\lambda\in\Lambda$ on the inside of $C_j(r)$. Also, on $C_{j-1}(r)$
we have
\[
U_j(z)\ge \frac{(\alpha+\varepsilon)\log\frac{1}{\eta(r)}+2C}
{\log\frac{1}{\delta(r)}+
        (\alpha+\varepsilon)\log\frac{1}{\eta(r)}}.
\]
It follows that
\[
R_j\ge \frac{\log\frac{1}{\delta(r)}-2C}{\log\frac{1}{\delta(r)}+
        (\alpha+\varepsilon)\log\frac{1}{\eta(r)}}.
\]
Thus
\[
\begin{split}
\log\frac{1}{P_n} &\le n\log \frac{\log\frac{1}{\delta(r)}+
        (\alpha+\varepsilon)\log\frac{1}{\eta(r)}}{\log\frac{1}{\delta(r)}-2C}\\
&= n \frac{(\alpha+\varepsilon)\log\frac{1}{\eta(r)}+2C}
          {\log\frac{1}{\delta(r)}-2C}
(1+o(1)) = (\alpha+\varepsilon) \frac{\log\frac{1}{1-r}}
                                     {\log\frac{1}{\delta(r)}}(1+o(1)).
\end{split}
\]
We finally have to estimate the infimum of the probabilities that
a particle starting from $C_n(r)$ hits $\partial \D$. Then take
$B$ to be the Blaschke product with zeros $z_l$, $|z_l|<r$, and
set
\[
U(z)=\frac{\log\frac{1}{|B(z)|}}{\log\frac{1}{\delta(r)}}
\]
Then since
\[
U(z)\le \frac{C}{\log\frac{1}{\delta(r)}} \quad \text{on } C_n(r),
\]
we get
\[
\log \frac{1}{\omega(\lambda,\partial \D,\Omega(\lambda,r))}
\le (\alpha+\varepsilon) \frac{\log\frac{1}{1-r}}{\log\frac{1}{\delta(r)}}
(1+o(1)).
\]

We now prove $\alpha \le D^+_h(\Lambda)$. Let $B_\zeta$ denote the
finite Blaschke product with zeros $\lambda\in\Lambda$ such that
$1/2<\rho(\zeta,\lambda)<r$. Define
\[
d(r)=\inf_{\lambda\in\Lambda} \log |B_{\lambda}(z)|
\]
and pick $\lambda^*$ such that
\[
\log|B_{\lambda^*}(\lambda^*)| < d(r) +1.
\]
By conformal invariance, we may assume $\lambda^*=0$.
We set
\[
d(r)+1=-(\alpha+\varepsilon)\log\frac{1}{1-r}
\]
and note that by our definition of $\alpha$,
$\varepsilon=\varepsilon(r)\to 0$ when $r\to 1$.

We introduce a function $\eta(r)$ as above and let $B_j$ and $U_j$
as before. We now build a harmonic function
\[
U=\sum_{j=1}^n w_n U_j
\]
such that $U(z)\le 1$ on $\partial D(\lambda,\delta(r))$.
First note that
\[
\sum_{j=1}^n U_j(z)\le
\log\frac{1}{\delta(z)} + (\alpha+\varepsilon)\log\frac{1}{\eta(r)}
\]
on the boundary of the ``bubbles" corresponding to the zeros of
$B_n$, with $\xi(r)\to 0$ as $\eta(r)\to 0$. Thus we set
\[
w_n=\frac{1}{a+b},
\]
where
\[
a= \log\frac{1}{\delta(r)}, \quad
b=(\alpha+\varepsilon)\log\frac{1}{\eta(r)}.
\]
Next we observe that on the boundary of the ``bubbles"
corresponding to the zeros of $B_{n-1}(z)$, we get
\[
w_{n-1}\sum_{j=1}^{n-1}U_{j}(z)+w_n U_n(z)\le w_{n-1}(a+b)
+ w_n b
\]
and so we set
\[
w_{n-1}=\frac{a}{(a+b)^2}.
\]
Inductively, we get
\[
w_{n-j}=\frac{1}{(a+b)}\left(\frac{a-b}a\right)^j.
\]
Note that the desired estimate on the boundaries of the ``bubbles"
is achieved because $w_j$ decreases when $j$ decreases.

To estimate $U(0)$, we argue in a similar way as above. The worst
case is that $U_n(0)$ is minimal because $w_n$ is the largest
weight. By our lower estimate
\[
\sum_{j=1}^n U_j(0)\le (\alpha+\varepsilon) \log\frac{1}{1-r}
\]
and the upper estimates $U_j(0)\le (\alpha+\xi)\log(1/\eta(r))$, we get
\[
U_j(0)\ge (\alpha+\xi)\log\frac{1}{\eta(r)}-
n(\xi-\varepsilon)\log\frac{1}{\eta(r)}.
\]
This leads us to the estimate
\[
U(0)\ge  1-\left(\frac{a}{a+b}\right)^n -
n(\xi+\varepsilon)\frac{\log\frac{1}{\eta(r)}}
{\log\frac{1}{\delta(r)}+\log\frac{1}{\eta(r)}},
\]
and so
\[
\omega(0,\partial \D,\Omega_\delta(0,r))\ge \left(\frac{a}{a+b}\right)^n
+ (\xi+\varepsilon)\frac{\log\frac{1}{1-r}}
{\log\frac{1}{\delta(r)}+\log\frac{1}{\eta(r)}}.
\]
We now require the second term to be ``small oh" of the first
term. This is certainly the case if $\delta(r)=1-r$. Thus
\[
\log\frac{1}{\omega(0,\partial \D,\Omega_\delta(0,r))}\ge
n \frac{b}{a}(1+o(1))=\alpha \frac{\log\frac{1}{1-r}}
{\log\frac{1}{\delta(r)}} (1+o(1)),
\]
and we are done.

\end{document}